\documentclass[10pt]{article}
\usepackage{graphicx}
\usepackage{subcaption}
\usepackage{bm}
\usepackage[authoryear]{natbib}
\usepackage{amsmath}
\usepackage[ruled,vlined]{algorithm2e}

\usepackage{indentfirst,csquotes}

\usepackage{amssymb,amsthm}
\usepackage{xcolor,paralist,hyperref,titlesec,fancyhdr,etoolbox}

\hypersetup{ colorlinks=true, linkcolor=black, filecolor=black, urlcolor=black }

\usepackage{lipsum}

\begin{document}
\title{Using Co-Located Range and Doppler Radars for Initial Orbit Determination} %%%%%%%%%%%%
\author{Cristina Parigini (1), Laura Pirovano (1), Roberto Armellin (1), \\ Darren McKnight (2), Adam Marsh (2), Tom Reddell (3), \\ 
((1) Te P\={u}naha \={A}tea - Space Institute, The University of Auckland, New Zealand, \\ (2) LeoLabs, United States, (3) LeoLabs, Australia)}  
\maketitle

%%%%%%%%%%

\begin{abstract}
With debris larger than 1 cm in size estimated to be over one million, precise cataloging efforts are essential to ensure space operations’ safety. Compounding this challenge is the oversubscribed problem, where the sheer volume of space objects surpasses ground-based observatories’ observational capacity. This results in sparse, brief observations and extended intervals before image acquisition. LeoLabs’ network of phased-array radars addresses this need by reliably tracking 10 cm objects and larger in low Earth orbit with 10 independent radars across six sites. While LeoLabs tracklets are extremely short, they hold much more information than typical radar observations. Furthermore, two tracklets are generally available, separated by a couple of minutes. Thus, this paper develops a tailored approach to initialize state and uncertainty from a single or pair of tracklets. Through differential algebra, the initial orbit determination provides the state space compatible with the available measurements, namely an orbit set. This practice, widely used in previous research, allows for efficient data association of different tracklets, thus enabling the addition of accurate tracks to the catalog following their independent initialization. The algorithm’s efficacy is tested using real measurements, evaluating the IOD solution’s accuracy and ability to predict the next passage from a single or a pair of tracklets.
\end{abstract} %%%%%%%%%

\bigskip

\noindent 
\section{Introduction}
Precise estimation of the state of Earth-orbiting objects is crucial for tasks such as observation scheduling, data association, and assessing collision risks. The reliability of these estimates largely depends on the sensor accuracy and the frequency of observations. Current resident space object catalogs only cover a small fraction of the total population, primarily objects larger than 10 cm. To enhance space safety, it is increasingly important to track smaller objects. This requires detection by sensors and the determination of their orbits with sufficient accuracy to ensure regular future observations and services such as conjunction predictions.

The task of estimating the first state of an uncatalogued object is known as initial orbit determination (IOD). Traditional IOD algorithms are classified based on the type of sensors and measurements involved. There are two primary categories: angle-only methods for optical sensors and angle-range methods for radar systems. The earliest solutions to the angle-only IOD problem were developed in the 18th century to calculate the orbits of celestial objects like planets and asteroids~\citep{Laplace1780, Gauss1857}. Recent improvements have been made to adapt to advancements in optical sensors~\citep{Escobal1965, BakerJr1977, Gooding1996, Karimi2010}. The second category of IOD methods is tailored for processing range-radar data, as discussed in~\citep{Lambert1761, Gibbs1889, Herrick1971, Zhang2020}, which use both range and angular measurements. In addition, specialized algorithms for handling Doppler and angular measurements have been proposed in recent years \citep{Losacco2023, Yanez2017}.

IOD methods generally do not provide uncertainty estimates unless a least squares refinement is applied or model simplifications are made~\citep{Zhang2020}. However, incorporating uncertainty quantification into the IOD solution is crucial for developing reliable data association strategies in catalog initialization. Recently, new approaches based on differential algebra (DA) have been proposed to overcome these shortcomings. Armellin and Di Lizia~\citep{Armellin2018} reformulated Lambert’s problem using DA, offering a mathematical representation of the uncertainty in the angle-range IOD solution as a function of measurement noise. This approach, when combined with automatic domain splitting (ADS)~\citep{Wittig2015}, was applied by \cite{Pirovano2020b} to solve the angle-only IOD problem and later by \citep{Losacco2023RobustIO} for Doppler-only radar systems.

Building on these works, this paper introduces IOD methods tailored to LeoLabs' phased-array radar network. LeoLabs’ sensors record range $\rho$, range-rate $\dot\rho$, azimuth $Az$, and elevation $El$ in low Earth orbit with 10 independent radars across six sites. Although the tracklets are very short, they contain more information than typical radar observations, necessitating a tailored approach to fully exploit the available data. This avoids using the classical admissible region approach for too-short arc, which would result in a larger uncertainty set. In most cases, each radar observation consists of two tracklets, each recorded by one of the two co-located arrays, separated by approximately 2 minutes. Furthermore, some sensors can obtain good-quality angular measurements for one of the two tracklets. Exploiting the geometry of the observations and the availability of ranges, range-rates and angle measurements, we hearby propose three algorithms.\\
\textit{Algorithm 1. Inaccurate angular measurements on one tracklet.}
In this first case, angular measurements are used only as a first guess to formulate a Lambert problem. The IOD problem is then solved by minimizing the residuals on ranges and range rates on the two tracks. \\
\textit{Algorithm 2. Accurate angular measurements on one tracklet.} On the other hand, when accurate angular measurements are provided, the IOD algorithm utilizes one set of ${\rho, \dot\rho, Az, El}$ and one set of $\rho, \dot\rho$ separated by 2 minutes to compute the orbit. \\
\textit{Algorithm 3. Angular measurements on both tracklets.} Lastly, assuming that both tracklets will have accurate angular measurements in the future, we propose a classical IOD method based on a Lambert problem across the two tracklets. The additional information available from the range-rates is used to reduce the orbital uncertainty. \\
Once a nominal IOD solution is obtained, DA techniques coupled with ADS are used to represent the IOD solution uncertainty and propagate it forward in time to the next observation opportunity, to enable association and then cataloging.

The paper is structured as follows: Section~\ref{Maths} describes DA and ADS. Section~\ref{Algo} provides the description of the three IOD algorithms. Finally, Section~\ref{sec:testcases} illustrates the algorithms' application to two objects in low Earth orbit.

\section{Mathematical tools}
\label{Maths}

\subsection{Differential Algebra}
\label{sec::DA}
Differential algebra is a mathematical framework that facilitates the computation of function derivatives within a computational setting. By replacing the algebra of real numbers with that of Taylor polynomials, any sufficiently regular function $\bm{f}$ of $v$ variables can be expanded into its Taylor series up to a desired order $k$. This algebra supports basic operations like addition, multiplication, as well as differentiation and integration~\citep{Berz1999}. A detailed explanation of DA applied to astrodynamics can be found in~\cite{DiLizia2008}. For instance, consider a multivariate function $\bm{y}=\bm{f}(\bm{x})$, with $\bm{x}\in\mathbb{R}^v$. Given a nominal value $\bm{x}$ and an associated uncertainty $\beta$, assumed uniform across all components, the DA representation of $\bm{x}$ can be expressed as
\begin{equation}
\label{eq:x_DA}
[\bm{x}]=\bm{x}+\beta\delta\bm{x}
\end{equation}
where $\delta\bm{x}\in[-1,1]^v$ represents the deviation from $\bm{x}$. When $\bm{f}$ is evaluated using the DA framework, the outcome is
\begin{equation}
\label{eq:y_DA}
[\bm{y}] = \bm{f}([\bm{x}])=\mathcal{T}_{\bm{y}}(\delta\bm{x})
\end{equation}
where $\mathcal{T}_{\bm{y}}(\delta\bm{x})$ in Eq.~\eqref{eq:y_DA} represents the Taylor expansion of $\bm{y}$ with respect to $\delta\bm{x}$. Therefore, for any initial deviation $\delta\bm{x}^{\ast}$, the corresponding result $\bm{y}^{\ast}$ can be determined by evaluating $\delta\bm{x}^{\ast}$ within the polynomial $\mathcal{T}_{\bm{y}}(\delta\bm{x})$, expressed as
\begin{equation}
\label{eq:y_star}
\bm{y}^{\ast}=\mathcal{T}_{\bm{y}}(\delta\bm{x}^{\ast})
\end{equation}

By combining DA with polynomial bounding techniques~\citep{Crane1975}, this approach offers a powerful means to estimate the bounds of $[\bm{y}]$.

\subsection{The ADS Algorithm}
\label{subsec:ADS}
The accuracy of the DA expansion generally decreases as the size of the uncertainty set $\delta\bm{x}$ grows or as the nonlinearity of $\bm{f}$ increases. One potential solution is to raise the expansion order $k$, but this doesn’t always lead to improved accuracy and often results in higher computational demands, making the DA approach impractical for highly nonlinear functions. However, the expansion order is not the only parameter for enhancing accuracy. An alternative approach involves subdividing the uncertainty set into smaller subsets, reducing the domain size. This technique recalculates the Taylor expansion around the center of each subset, maintaining overall coverage but achieving higher accuracy within each expansion. Following these principles, ADS uses an automatic algorithm to assess whether the current polynomial representation meets the required accuracy~\citep{Wittig2015}. If the desired accuracy is not reached, the original domain is split in half, and two new expansions are generated by re-expanding the polynomials around the new centers, each covering half of the initial set. This process is repeated for the newly created subsets until the target accuracy is achieved or a predefined limit on the number of splits is reached. The result is a \textit{manifold} of Taylor expansions, represented as
\begin{equation}
\bm{y}=\bigcup\limits_{i=1}^{N_s}\mathcal{T}_{\bm{y}}^i(\delta\bm{x}).
\end{equation}

\subsection{LeoLabs' radar observations}
\label{sec::LR}
When an object is detected at one of the LeoLabs sites, it typically produces two tracklets by passing through two co-located radars. Each tracklet comprises a set of observations within a short time span of a few seconds, separated by a couple of minutes. While an individual tracklet can be used to determine the immediate passage to the co-located radar, the resulting orbital uncertainty would be too large for data association at the next observation window, typically separated by a few hours. 
\begin{figure}[!b]
    \centering
    \includegraphics[width=\columnwidth]{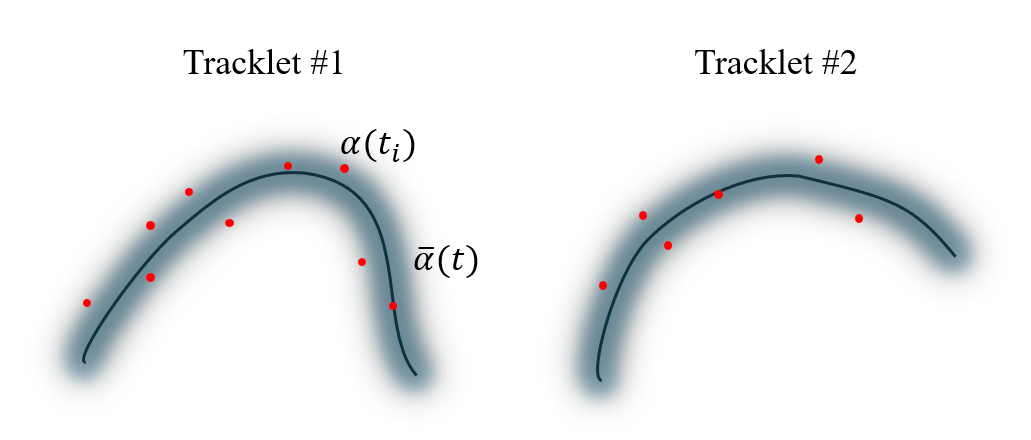}  % Replace with your image path
    \caption{Illustration of the two tracklets.}
    \label{fig:tracklet}  % Replace with a specific label for referencing
\end{figure}
Referring to Fig. \ref{fig:tracklet} a generic measurement part of the first or the second tracklet is indicated with $\alpha(t_i)$, with $\alpha \in \{\rho, \dot\rho, Az, El\}$. These are raw measurements taken at the time $t_i$ (red dots). To account for sensor level errors, the observation's precision can be modeled as white noise and thus be considered a Gaussian random variable with zero mean and \(\sigma\) standard deviation \citep{Numerica}. 

A estimate of the measurements obtained by a polynomial regression of all the data in a tracklet is indicated with $\bar{\alpha}(t_i)$ (for details about the regression, see \cite{Losacco2023} or \cite{Pirovano2020b}). 
When the regressed data are used, their precision is estimated by the confidence interval (CI) of the regression for a chosen significance level $\alpha$. Regression and CI are sketched as the continuous curve and its shade in Fig. \ref{fig:tracklet}.
Through the time derivative of regression polynomials, quantities not directly measured, i.e. $\dot{\bar{Az}},\dot{\bar{El}}$, can be estimated. However, due to the relatively low accuracy in angular measurements, these quantities, when used, can provide only rough initial guesses for the IOD algorithms.  

\section{IOD Algorithms} \label{Algo}
In this section, we describe three main algorithms that were developed and tested based on the different observations available in the two tracklets. It is worth mentioning that the first option that we considered was to derive azimuth and elevation rates from interpolated azimuth and elevation data, similar to what was proposed in \cite{d2023initial}, but the results we obtained were largely unreliable, probably due to the shortness of the tracklets and the accuracy of the measurements.  

\subsection{Inaccurate angular measurements on one tracklet}
In case the angular measurements are inaccurate, we can only rely on pairs of range and range-rate measurements with three of them, at time $t_1$, $t_2$ and $t_3$ belonging to one tracklet. If measurements from the other tracklet are available, a fourth set of range-range rate data is considered at time $t_4$. The IOD algorithm is summarized in Algorithm \ref{alg:IOD1}.

\begin{algorithm}[H]
\caption{IOD with inaccurate angular measurements}\label{alg:IOD1}
\KwIn{Initial guesses for $\bar{Az}(t_1)$, $\bar{El}(t_1)$, $\bar{Az}(t_3)$, $\bar{El}(t_3)$}
\KwOut{The values of $\bar{Az}(t_1)$, $\bar{El}(t_1)$, $\bar{Az}(t_3)$, $\bar{El}(t_3)$ that bring the residuals to zero $\boldsymbol{\Delta} = \boldsymbol{0}$ or that minimize $J = \boldsymbol{\Delta}^T \boldsymbol{\Delta}$}

\BlankLine

1. Guess initial values for $\bar{Az}(t_1)$, $\bar{El}(t_1)$, $\bar{Az}(t_3)$, $\bar{El}(t_3)$\;

2. Take observations $\rho(t_1)$ and $\rho(t_3)$\;

3. Determine the initial position vector: \\
$\boldsymbol{r}(t_1)(\rho(t_1), \bar{Az}(t_1), \bar{El}(t_1))$\;

4. Determine the final position vector: \\
$\boldsymbol{r}(t_3)(\rho(t_3), \bar{Az}(t_3), \bar{El}(t_3))$

5. Solve Lambert's problem to compute: \\
$\boldsymbol{v}(t_1)(\rho(t_1), \bar{Az}(t_1), \bar{El}(t_1), \rho(t_3), \bar{Az}(t_3), \bar{El}(t_3))$ \\
$\boldsymbol{v}(t_3)(\rho(t_1), \bar{Az}(t_1), \bar{El}(t_1), \rho(t_3), \bar{Az}(t_3), \bar{El}(t_3))$ \;

6. Propagate the state $(\boldsymbol{r}(t_1), \boldsymbol{v}(t_1))$ to an intermediate time $t_2$ using Keplerian dynamics, and, if data from a second tracklet are available to $t_4$: \\
$\boldsymbol{r}(t_2)(\rho(t_1), \bar{Az}(t_1), \bar{El}(t_1), \rho(t_3), \bar{Az}(t_3), \bar{El}(t_3))$ \\
$\boldsymbol{v}(t_2)(\rho(t_1), \bar{Az}(t_1), \bar{El}(t_1), \rho(t_3), \bar{Az}(t_3), \bar{El}(t_3))$\\
$\boldsymbol{r}(t_4)(\rho(t_1), \bar{Az}(t_1), \bar{El}(t_1), \rho(t_3), \bar{Az}(t_3), \bar{El}(t_3))$ \\
$\boldsymbol{v}(t_4)(\rho(t_1), \bar{Az}(t_1), \bar{El}(t_1), \rho(t_3), \bar{Az}(t_3), \bar{El}(t_3))$\;

7. Build the residual vector $\boldsymbol{\Delta}$ using the states at $t_1, t_2, t_3$ and, if available, at $t_4$: \\
$
\boldsymbol{\Delta}(\rho(t_1), \bar{Az}(t_1), \bar{El}(t_1), \rho(t_3), \bar{Az}(t_3), \bar{El}(t_3)) =
\begin{pmatrix}
\rho^c(t_1) - \rho(t_1) \\
\dot{\rho}^c(t_1) - \dot{\rho}(t_1) \\
\rho^c(t_2) - \rho(t_2) \\
\dot{\rho}^c(t_2) - \dot{\rho}(t_2) \\
\rho^c(t_3) - \rho(t_3) \\
\dot{\rho}^c(t_3) - \dot{\rho}(t_3) \\
\rho^c(t_4) - \rho(t_4) \\
\dot{\rho}^c(t_4) - \dot{\rho}(t_4) \\
\end{pmatrix}$\; 
where the superscript $c$ indicates computed quantities to differentiate them from the measurements. 

8. Compute the values $\bar{Az}(t_1)$, $\bar{El}(t_1)$, $\bar{Az}(t_3)$, $\bar{El}(t_3)$ such that $\boldsymbol{\Delta} = \boldsymbol{0}$ or to minimize $J = \boldsymbol{\Delta}^T \boldsymbol{\Delta}$

\end{algorithm}

The solution to the set of nonlinear equations or the minimization problem is obtained with the Broyden–Fletcher–Goldfarb–Shanno algorithm implemented in dlib library \cite{king2009dlib}. However, while this algorithm can provide an IOD solution, in practice it is not useful as it fails to provide the entire set of orbits compatible with the measurements and their uncertainty, i.e. the orbit set. To solve this issue Algorithm \ref{alg:IOD1} is run in the DA framework using ADS. We start by initializing the variables as 

\begin{equation}
\begin{aligned}
&[\rho(t_1)] = \rho(t_1) &+ 3\sigma_{\rho(t_1)} \cdot \delta \rho(t_1) \\
&[\bar{Az}(t_1)] = \bar{Az}(t_1) &+ CI_{\bar{Az}(t_1)} \cdot \delta \bar{Az}(t_1)\\
&[\bar{El}(t_1)] = \bar{El}(t_1) &+ CI_{\bar{El}(t_1)} \cdot \delta \bar{El}(t_1)\\
&[\rho(t_3)] = \rho(t_3) &+ 3\sigma_{\rho(t_3)} \cdot \delta \rho(t_3) \\
&[\bar{Az}(t_3)] = \bar{Az}(t_3) &+ CI_{\bar{Az}(t_3)} \cdot \delta \bar{Az}(t_3)\\
&[\bar{El}(t_3)] = \bar{El}(t_3) &+ CI_{\bar{El}(t_3)} \cdot \delta \bar{El}(t_3)\\
\end{aligned}
\end{equation}
where $CI$ indicates the confidence interval for the measurements or search space for the optimization variables. We then perform steps 1--7 of Algorithm \ref{alg:IOD1} with ADS. This operation  produces  
\begin{equation}
\bm{\Delta}^T\bm{\Delta}=\bigcup\limits_{i=1}^{N_s}\mathcal{T}_{\bm{\Delta}^T \bm{\Delta}}^i(\bm{x})
\end{equation}
with $ \bm{x} = (\rho(t_1), \bar{Az}(t_1), \bar{El}(t_1), \rho(t_3), \bar{Az}(t_3), \bar{El}(t_3))$ and $N_s$ the number of domains generated by ADS. 
A minimization of the residual squares is then performed on each domain. Only those for which $\bm{\Delta}^T\bm{\Delta}$ is zero or below a certain thresholds are retained. 
These domains describe the set of spacecraft states compatible with the available observations, which can propagated forward in time to predict future detections. 
  
\subsection{Accurate angular measurements on one tracklet}
The second case considers the availability of accurate angular measurements on one of the two tracklets. For illustration, we consider these measurements together with the range and range-rate to be available on the first tracklet at $t_1$. Range and range-rate are also available for the second trackelet, at $t_2$. Concerning the angular data and their rates, we work with the regressed values and their rate.  The procedure is summarized in Algorithm \ref{alg:IOD2}.

\begin{algorithm}[ht]
\caption{IOD with accurate angular measurements on a single tracklet}
\label{alg:IOD2}
\KwIn{Initial guesses for $\dot{\bar{Az}}(t_1)$, $\dot{\bar{El}}(t_1)$}
\KwOut{The values of $\dot{\bar{Az}}(t_1)$, $\dot{\bar{El}}(t_1)$ that solve $\boldsymbol{\Delta} = \boldsymbol{0}$}

\BlankLine

1. Guess initial values for $\dot{\bar{Az}}(t_1)$, $\dot{\bar{El}}(t_1)$\;

2. Take observations ${{\rho}}(t_1)$, $\dot{{\rho}}(t_1)$, ${\bar{Az}}(t_1)$, ${\bar{El}}(t_1)$\;

3. Compute the Cartesian state: \\
   $\boldsymbol{r}(t_1)({{\rho}}(t_1), {{\rho}}(t_1), {\bar{Az}}(t_1), {\bar{El}}(t_1), \dot{\bar{Az}}(t_1), \dot{\bar{El}}(t_1))$\\
   $\boldsymbol{v}(t_1)({{\rho}}(t_1), \dot{{\rho}}(t_1), {\bar{Az}}(t_1), {\bar{El}}(t_1), \dot{\bar{Az}}(t_1), \dot{\bar{El}}(t_1))$\;

4. Propagate to the the state forward in time: \\
    $\boldsymbol{r}(t_2)({{\rho}}(t_1), \dot{{\rho}}(t_1), {\bar{Az}}(t_1), {\bar{El}}(t_1), \dot{\bar{Az}}(t_1), \dot{\bar{El}}(t_1))$\\
    $\boldsymbol{v}(t_2)({{\rho}}(t_1), \dot{{\rho}}(t_1), {\bar{Az}}(t_1), {\bar{El}}(t_1), \dot{\bar{Az}}(t_1), \dot{\bar{El}}(t_1))$\;

5. Build the residual vector $\boldsymbol{\Delta}$ using the state at $ t_2$: \\
$\boldsymbol{\Delta}({{\rho}}(t_1), \dot{{\rho}}(t_1), {\bar{Az}}(t_1), {\bar{El}}(t_1), \dot{\bar{Az}}(t_1, \dot{\bar{El}}(t_1))) =    \begin{pmatrix}
   \rho^c(t_2) - {\rho}(t_2) \\
   \dot{\rho}^c(t_2) - \dot{{\rho}}(t_2) 
   \end{pmatrix}$\;

6. Compute $\dot{\bar{Az}}(t_1)$, $\dot{\bar{El}}(t_1)$ such that $\boldsymbol{\Delta} = \boldsymbol{0}$ 
\end{algorithm}

Algorithm \ref{alg:IOD2} solves the IOD as a 2D nonlinear problem. If uncertainties in the observation are taken into account, then, by using the high-order Taylor expansion of the solution of the residual constraint, it is possible to compute the orbit set \cite{Pirovano2020b}. To do so, we first initialize   

\begin{equation}
\begin{aligned}
\label{eq:ini}
&[{\rho}(t_1)] = {\rho}(t_1) &+ 3 \sigma_{{\rho}(t_1)} \cdot \delta {\rho}(t_1) \\
&[\bar{Az}(t_1)] = \bar{Az}(t_1) &+ CI_{\bar{Az}(t_1)} \cdot \delta \bar{Az}(t_1)\\
&[\bar{El}(t_1)] = \bar{El}(t_1) &+ CI_{\bar{El}(t_1)} \cdot \delta \bar{El}(t_1) \\
&[\dot{{\rho}}(t_1)] = \dot{{\rho}}(t_1) &+ 3 \sigma_{\dot{\bar{\rho}}(t_1)} \cdot \delta \dot{{\rho}}(t_1) \\
&[{\rho}(t_2)] = {\rho}(t_2) &+ 3 \sigma_{{\rho}(t_2)} \cdot \delta {\rho}(t_2) \\
&[\dot{{\rho}}(t_2)] = \dot{{\rho}}(t_2) &+ 3 \sigma_{\dot{\bar{\rho}}(t_2)} \cdot \delta \dot{{\rho}}(t_2) \\
\end{aligned}
\end{equation}
Assuming an initial guess for $\dot{\bar{Az}}(t_1)$ and $\dot{\bar{El}}(t_1)$ based on the regression, we compute the Taylor expansion of the computed range and range rate at the second tracklet,
\begin{equation}
\begin{aligned}
\rho^c(t_2) = \mathcal{T}_{\rho^c(t_2)} \Big( {{\rho}}(t_1), \dot{{\rho}}(t_1), {\bar{Az}}(t_1), {\bar{El}}(t_1), \dot{\bar{Az}}(t_1), \dot{\bar{El}}(t_1) \Big) \\
\dot{\rho}^c(t_2) = \mathcal{T}_{\dot{\rho}^c(t_2)} \Big( {{\rho}}(t_1), \dot{{\rho}}(t_1), {\bar{Az}}(t_1), {\bar{El}}(t_1), \dot{\bar{Az}}(t_1), \dot{\bar{El}}(t_1) \Big)
\end{aligned}
\end{equation}

Partial inversion techniques are then used to enforce the computed range and range rates to be equal to the actual measurements 

\begin{equation}
\begin{aligned}
\label{eq:az}
\dot{\bar{Az}}(t_1) &= \mathcal{T}_{\dot{\bar{Az}}(t_1)} \Big( \rho^c(t_2) = {\rho}(t_2), \dot{\rho}^c(t_2) = \dot{{\rho}}(t_2), \\
&\quad {{\rho}}(t_1), \dot{{\rho}}(t_1), {\bar{Az}}(t_1), {\bar{El}}(t_1) \Big) \\
&= \mathcal{T}_{\dot{\bar{Az}}(t_1)} \Big( {{\rho}}(t_1), \dot{{\rho}}(t_1), {\bar{Az}}(t_1), {\bar{El}}(t_1), \\
&\quad {\rho}(t_2), \dot{{\rho}}(t_2) \Big).
\end{aligned}
\end{equation}

\begin{equation}
\begin{aligned}
\label{eq:el}
\dot{\bar{El}}(t_1) &= \mathcal{T}_{\dot{\bar{El}}(t_1)} \Big( \rho^c(t_2) = {\rho}(t_2), \dot{\rho}^c(t_2) = \dot{{\rho}}(t_2), \\
&\quad {{\rho}}(t_1), \dot{{\rho}}(t_1), {\bar{Az}}(t_1), {\bar{El}}(t_1) \Big) \\
&= \mathcal{T}_{\dot{\bar{El}}(t_1)} \Big( {{\rho}}(t_1), \dot{{\rho}}(t_1), {\bar{Az}}(t_1), \\
&\quad {\bar{El}}(t_1), {\rho}(t_2), \dot{{\rho}}(t_2) \Big).
\end{aligned}
\end{equation}

Equations \eqref{eq:az} and \eqref{eq:el}, along with the first four in Eq. \ref{eq:ini}, define the set of orbits that are consistent with the observations ${{\rho}}(t_1)$, ${\dot{\rho}}(t_1)$, ${\bar{Az}}(t_1)$, ${\bar{El}}(t_1)$, ${\rho}(t_2)$, $\dot{{\rho}}(t_2)$ and their uncertainties. By conducting all calculations within the ADS, the orbit set is represented as a series of Taylor polynomials with specified accuracy. This polynomial representation can be propagated forward in time and projected into the observation space when additional observations become available, aiding in data association and the cataloging of new objects.

\subsection{Accurate angular measurements on two tracklets}
When accurate angular measurements are available on both tracklets, the IOD problem is firstly tackled as a standard range radar, solved via the expansion of the Lambert problem solution with respect to the measurements  ${{\rho}}(t_1), {\bar{Az}}(t_1), {\bar{El}}(t_1), {\rho}(t_2), {\bar{Az}}(t_2), {\bar{El}}(t_2)$, as in \cite{Armellin2018}, in which $t_1$ and $t_2$ are respectively in the first and second tracklet. To exploit the additional accurate information provided by the range rates measurement, the IOD solver is evaluated in the ADS framework, and computed range rates are derived from the velocities provided by the Lambert problem as 

\begin{equation}
\dot{\rho}^c(t_1) = \bigcup\limits_{i=1}^{N_s} \mathcal{T}_{\dot{\rho}^c(t_1)}^i(\bm{x})
\end{equation}
\begin{equation}
\dot{\rho}^c(t_2) = \bigcup\limits_{i=1}^{N_s} \mathcal{T}_{\dot{\rho}^c(t_2)}^i(\bm{x})
\end{equation}
with $\boldsymbol{x} = ({{\rho}}(t_1), {\bar{Az}}(t_1), {\bar{El}}(t_1), 
{\rho}(t_2), {\bar{Az}}(t_2), {\bar{El}}(t_2))$.
Domains in which the computed range rates are incompatible with the measurements are pruned away, thus drastically reducing the size of the orbit set. We refer to this method as Algorithm 3 in the reminder of the paper. 

\section{Test cases}\label{sec:testcases}
\begin{figure*}[]
     \centering
     \begin{subfigure}[b]{0.48\textwidth}
         \centering
        \includegraphics[width=\textwidth]{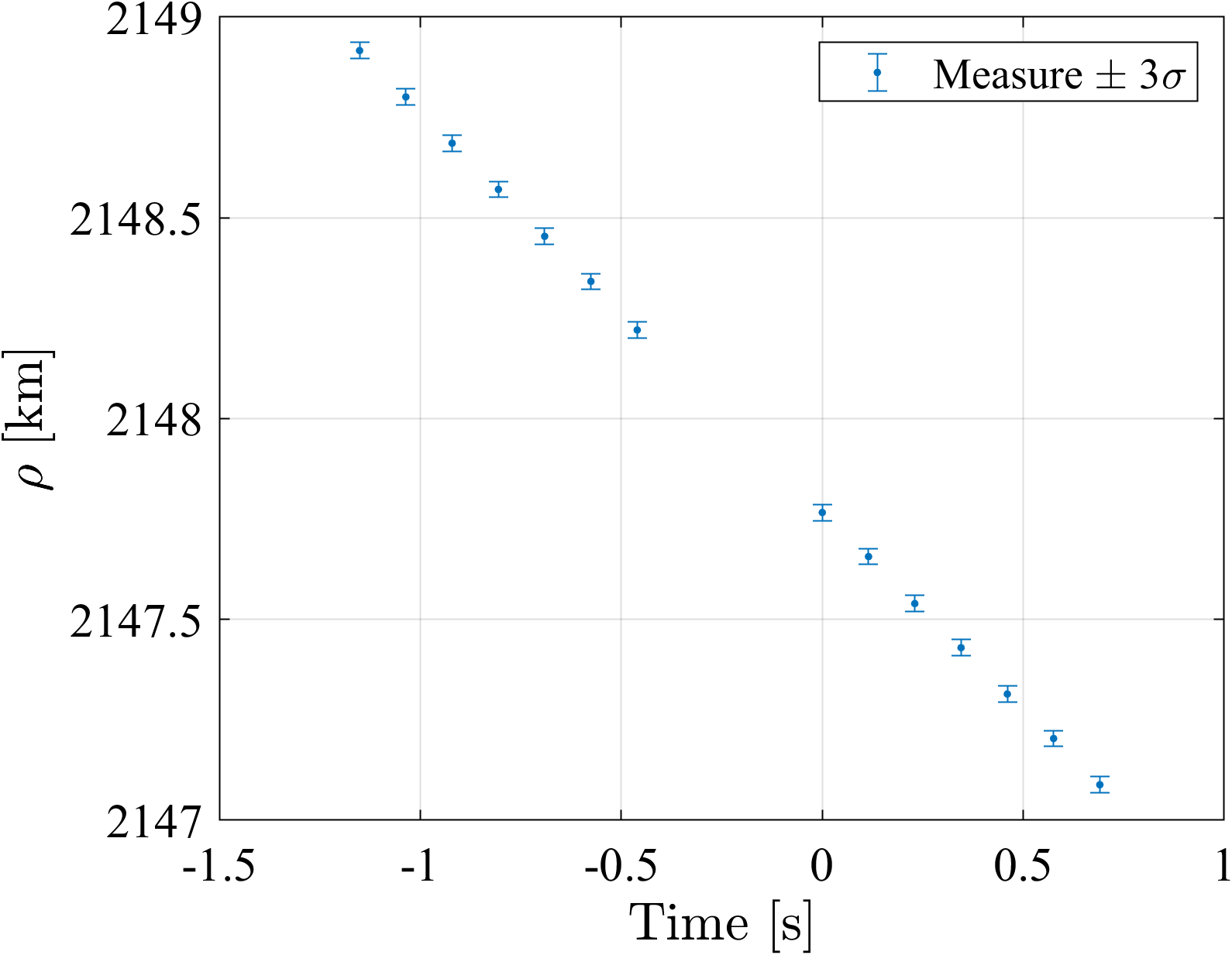}
         \caption{Range measurements and uncertainty.}
         \label{fig:data_range}
     \end{subfigure}
     \hfill
     \begin{subfigure}[b]{0.48\textwidth}
         \centering
        \includegraphics[width=\textwidth]{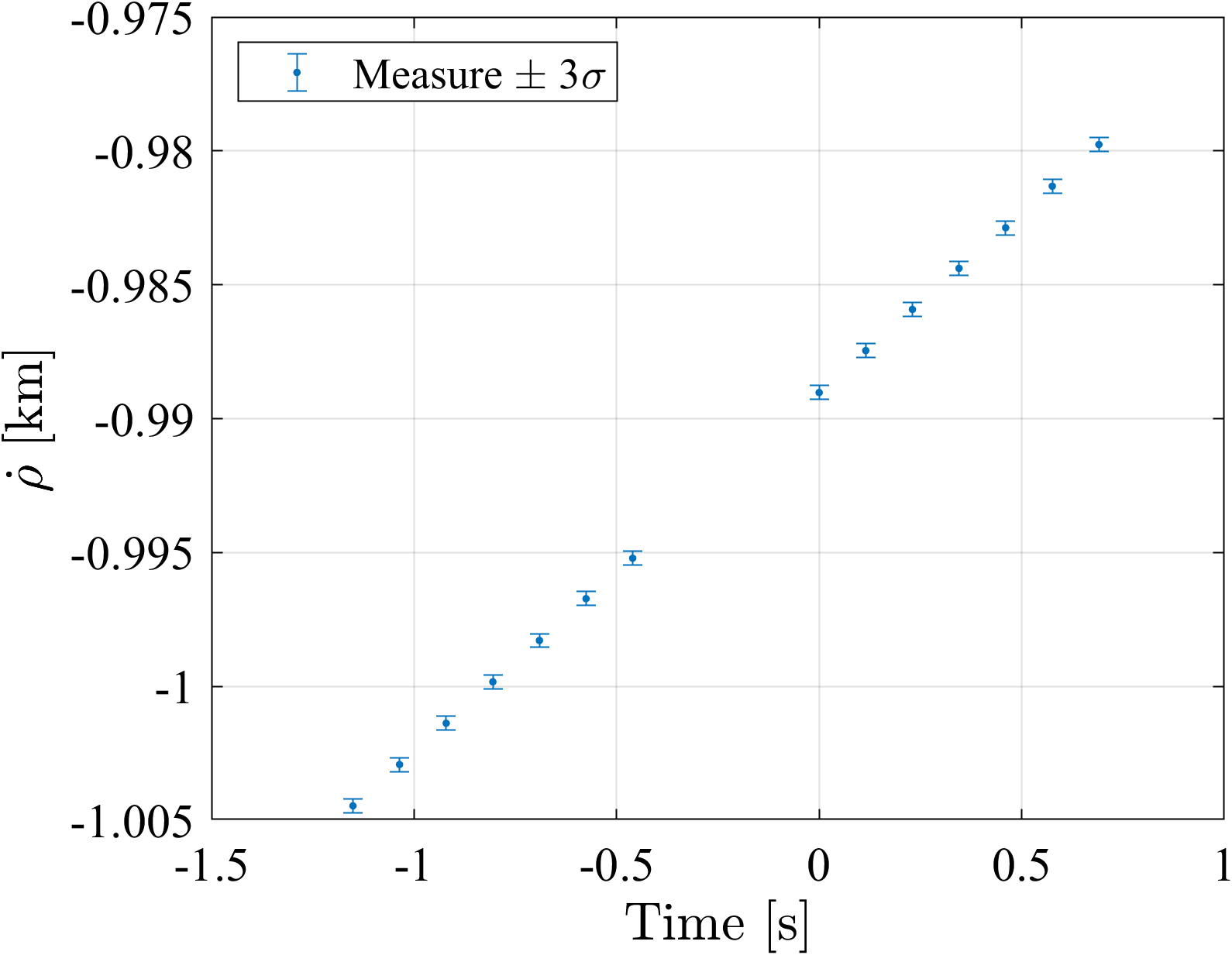}
         \caption{Range rate measurements and uncertainty.}
         \label{fig:data_rangeRate}
     \end{subfigure}
     \newline
     \begin{subfigure}[b]{0.48\textwidth}
         \centering
         \includegraphics[width=\textwidth]{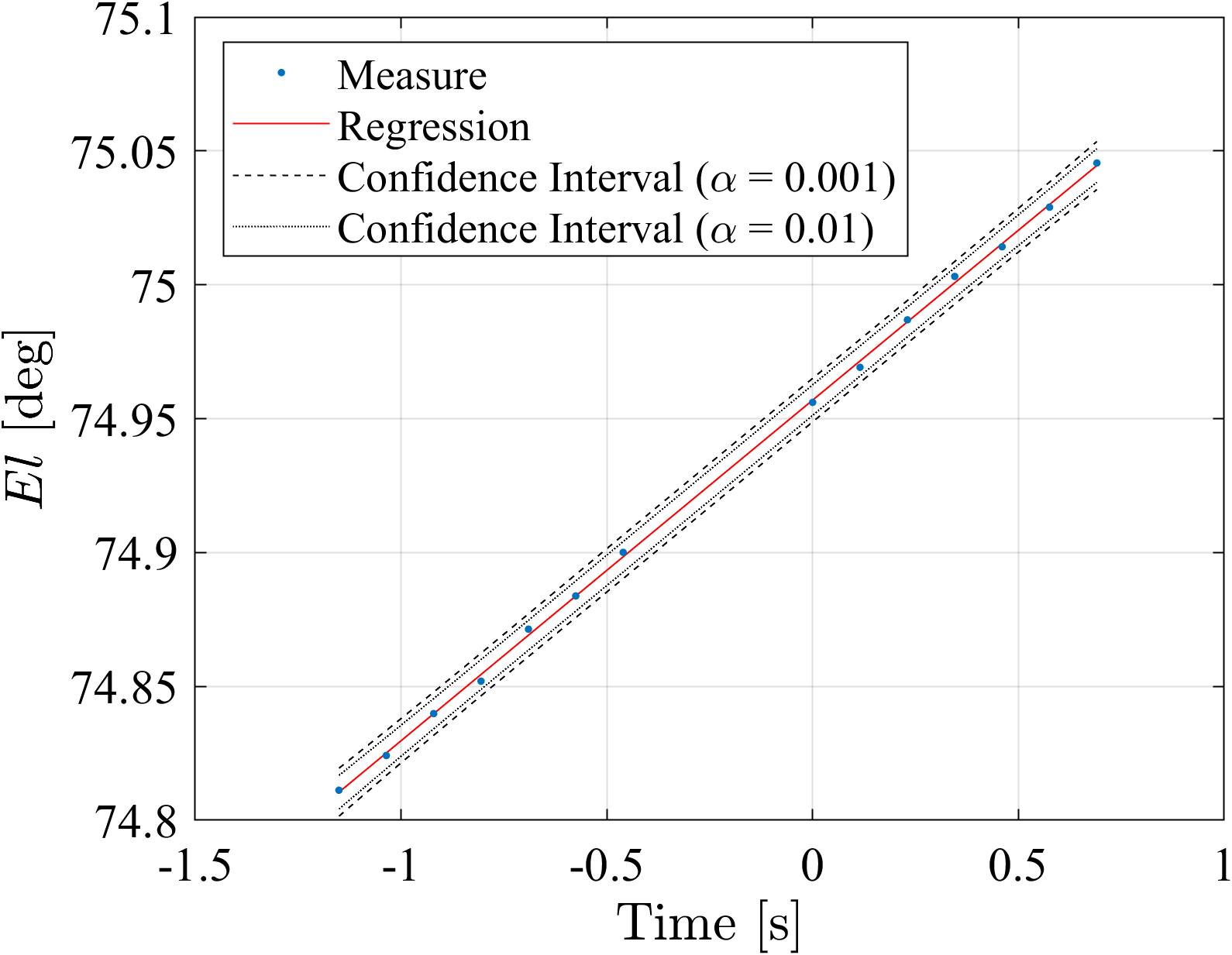}
         \caption{Azimuth measurements and regression.}
         \label{fig:data_azimuth}
     \end{subfigure}
     \hfill
     \begin{subfigure}[b]{0.48\textwidth}
         \centering
         \includegraphics[width=\textwidth]{figures/data_L12223_T2_elevation.png}
    \caption{Elevation measurements and regression.}
    \label{fig:data_elevation}
     \end{subfigure}
        \caption{Raw data and polynomial regression on Object 1, tracklet 2.}
        \label{fig:Reg1}
\end{figure*}

 Two objects with NORAD IDs 10011 and 39027 observed by the Costa Rica Space Radar are used as test cases. The orbital parameters for these objects are detailed in Table \ref{tab:orbital_params}, along with the epoch at the end of their second tracklet. Each object has three associated tracklets: the first lacks angular information, the second includes angular information and is separated from the first by roughly two minutes, indicating they are part of the same passage, and the third occurs around 12 hours later, marking the subsequent radar passage. The third tracklet's purpose is to evaluate the feasibility of data association, while the first two tracklets are utilized for IOD.
 Data for tracklet 2 of Object 1 is shown in Fig.~\ref{fig:Reg1}. Range and range rate measurements, along with their 3$\sigma$ values, are directly provided by LeoLabs. For angular measurements, due to the absence of provided uncertainties, a polynomial regression is performed to determine the CI, following the methodology in \cite{Pirovano2020b}. Two significance levels, $\alpha$, are used to obtain the CI. In all cases in this study, it is presumed that there is no prior knowledge of the objects' orbital parameters. LeoLabs' accurate catalog data is solely for validation. 

\begin{table}[]
    \centering
    \begin{tabular}{|c|c|c|}
        \hline
        Parameter & Object 1 & Object 2 \\
        \hline
        NORAD & 10011 & 39027 \\
        Date & 2024-05-03 & 2024-04-29 \\
        Time (UTC) & 12:15:07.052812Z & 09:47:30.383527Z \\
        $a$ (km) & 7921.54 & 6849.12 \\
        $e$ & 0.06974 & 0.00553 \\
        $i$ (deg) & 65.9450 & 97.6549 \\
        $\Omega$ (deg) & 147.613 & 98.4364 \\
        $\omega$ (deg) & 338.098 & 118.156 \\
        $\vartheta$ (deg) & 186.836 & 50.813 \\
        \hline
    \end{tabular}
    \caption{Orbital parameters for 10011 and 39027.}
    \label{tab:orbital_params}
\end{table}

% piccolo lambert
\subsection{Association of very short tracklets}
\begin{figure*}[h]
     \centering
     \begin{subfigure}[b]{0.48\textwidth}
         \centering
         \includegraphics[width=\columnwidth]{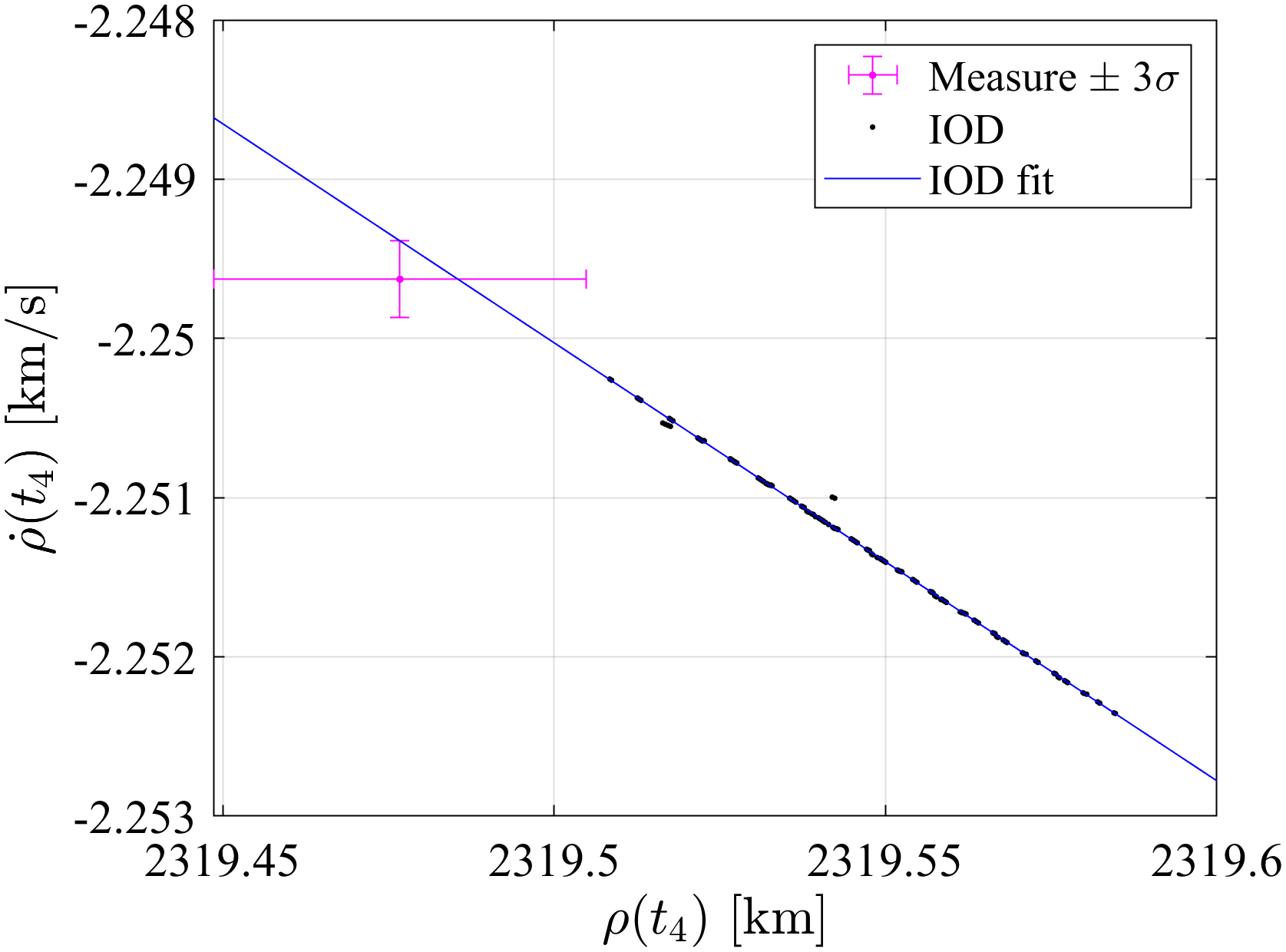}
    \caption{Object 1.}
    \label{fig:data_lambertinoObj1}
     \end{subfigure}
     \hfill
     \begin{subfigure}[b]{0.48\textwidth}
         \centering
        \includegraphics[width=\columnwidth]{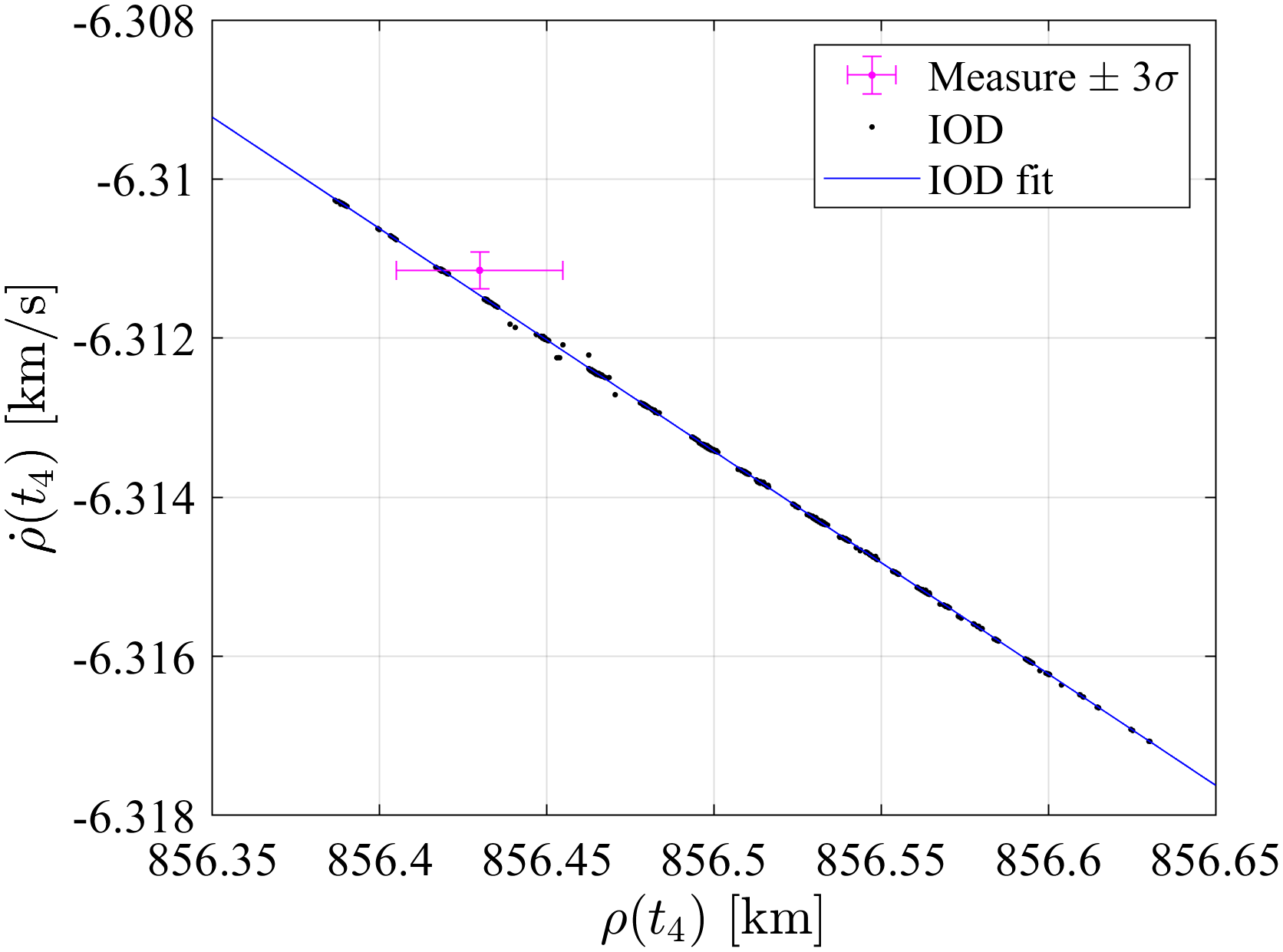}
    \caption{Object 2.}
    \label{fig:data_lambertinoObj2}
     \end{subfigure}
     \caption{Association analysis on tracklet 2 for the two cases.}
        \label{fig:Assoc}
\end{figure*}

The initial analysis is conducted to determine whether two tracklets, lasting a few seconds and separated by approximately two minutes, can be associated. For this analysis, we use Algorithm \ref{alg:IOD1}, considering only the range and range-rate measurements within each tracklet to construct the residual function $J$, as it is unknown at this stage whether the second tracklet belongs to the same object. The states corresponding to the minimum residuals, computed within each domain generated by the ADS, are propagated forward/backward to the time of the second tracklet and compared with new measurements. Figures \ref{fig:data_lambertinoObj1} and \ref{fig:data_lambertinoObj2} present the results for the two test cases. In both cases, it can be observed that the IOD solutions yield the lowest residuals along a line in the range and range-rate plane, consistent with the new observations taken during the radar’s second pass, shown in pink. This suggests a potential association between the tracklets, although this may not be the only possibility, depending on the density of uncorrelated tracks in the measurement plane.

\subsection{IOD with inaccurate angular measurements} % lambert
Based on the assumption that the first two tracklets are successfully associated, we now use Algorithm \ref{alg:IOD1} to perform the IOD using both tracklets. Algorithm \ref{alg:IOD1} is run to compute the solution, considering a significance level $\alpha = 0.001$ on the regressed $Az$ and $El$ data. This large confidence interval is used to account for inaccurate measurements. 

\begin{figure*}[h]
     \centering
     \begin{subfigure}[b]{0.48\textwidth}
         \centering
        \includegraphics[width=\columnwidth]{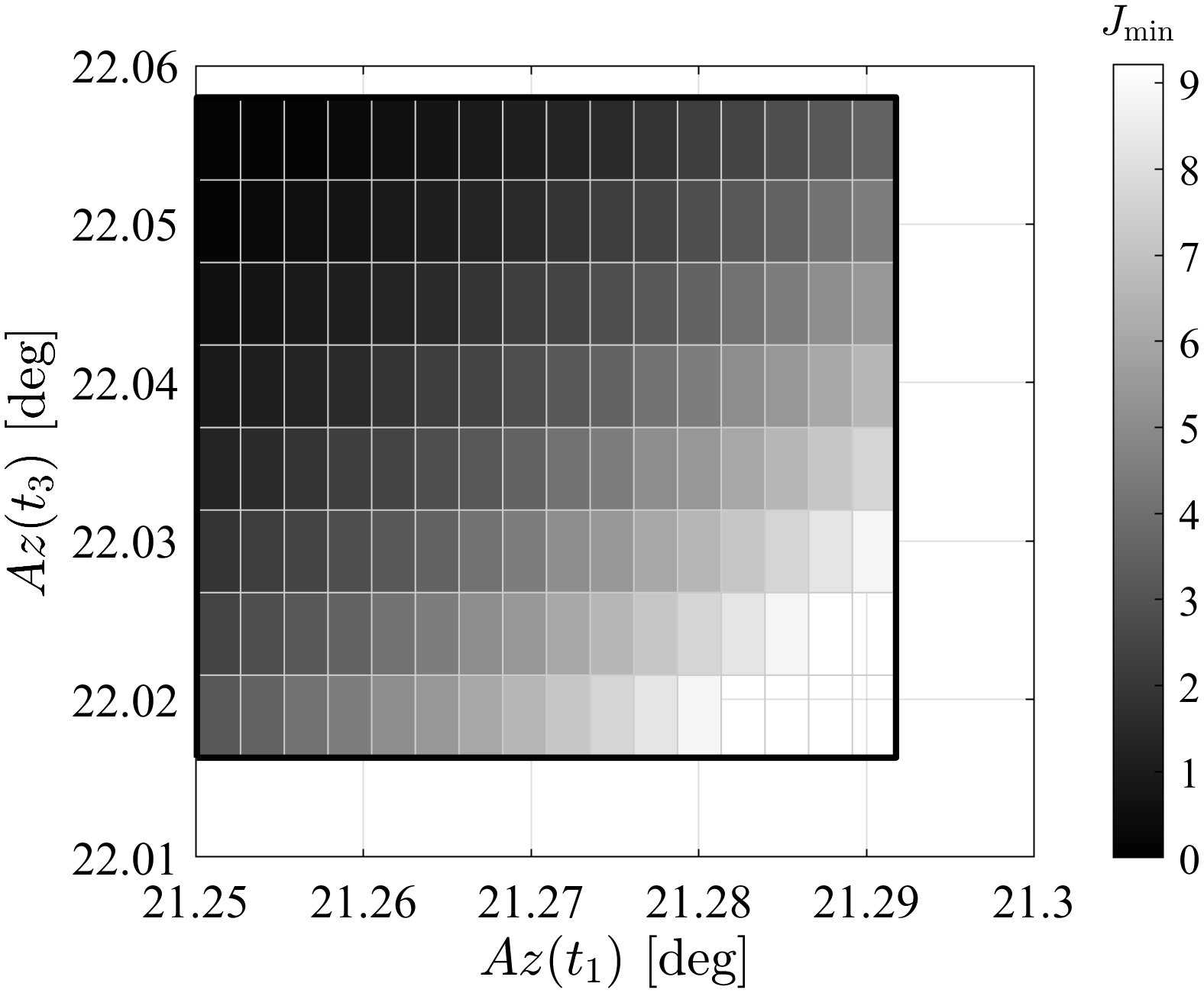}
    \caption{Domains in the $Az$ space.}
    \label{fig:LambertDomainA1A3_obs2m_TC1}
     \end{subfigure}
     \hfill
     \begin{subfigure}[b]{0.48\textwidth}
         \centering
        \includegraphics[width=\columnwidth]{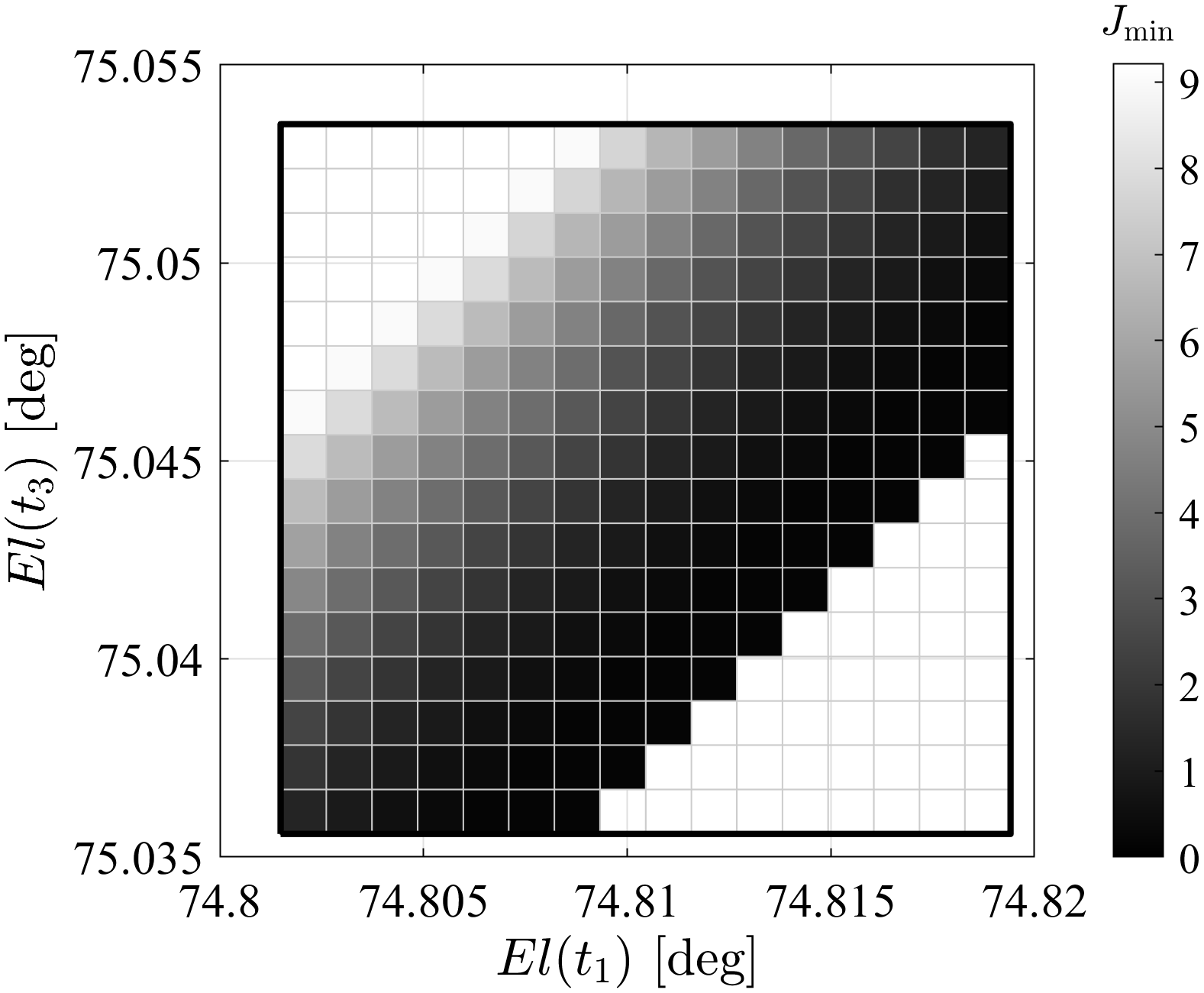}
    \caption{Domains in the $El$ space}
    \label{fig:LambertDomainE1E3_obs2m_TC1}
     \end{subfigure}
     \caption{Domain and residual for Object 1, Algorithm 1.}
        \label{fig:ResultAlgo1}
\end{figure*}

The results are reported only for Object 1, as both cases produce similar results. Figure \ref{fig:LambertDomainA1A3_obs2m_TC1} and \ref{fig:LambertDomainE1E3_obs2m_TC1} report the domains generated by ADS in the angular measurements space, colored by the minimum value reached by the target function $J_{min}$.  Note that the domains are defined in a 6D space; thus, multiple domains are present in other dimensions in the two figures. The color code is used to indicate the minimum across all the dimensions. White domains are those pruned away considering as condition $J < \chi^2_{0.99, 2} = 9.21$, where $\chi^2_{0.99, 2}$ is 99th percentile of the Chi-squared distribution with 2 degrees of freedom (as we are considering 8 measurements in total to build the residual function).
Using a tolerance of $10^{-5}$ and a maximum number of splits of $15$, in both cases, the maximum number of domains of $32,768$ is reached. The pruning method then leaves 5,403 feasible boxes for the first test case and 891 for the second.  

\begin{figure*}[h!]
     \centering
     \begin{subfigure}[b]{0.48\textwidth}
         \centering
       \includegraphics[width=\columnwidth]{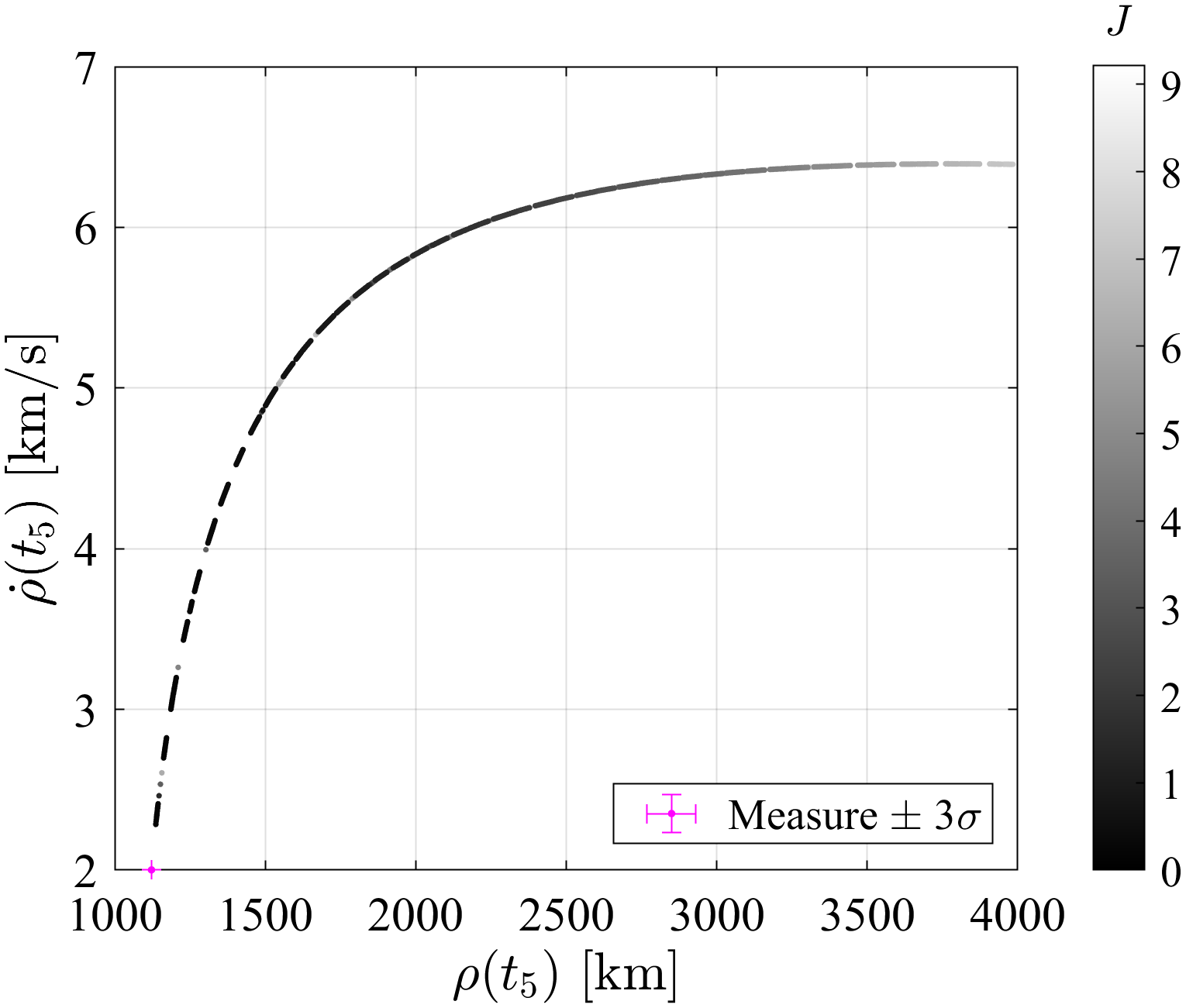}
    \caption{Object 1.}
    \label{fig:LambertPropag12_TC1}
     \end{subfigure}
     \hfill
     \begin{subfigure}[b]{0.48\textwidth}
         \centering
        \includegraphics[width=\columnwidth]{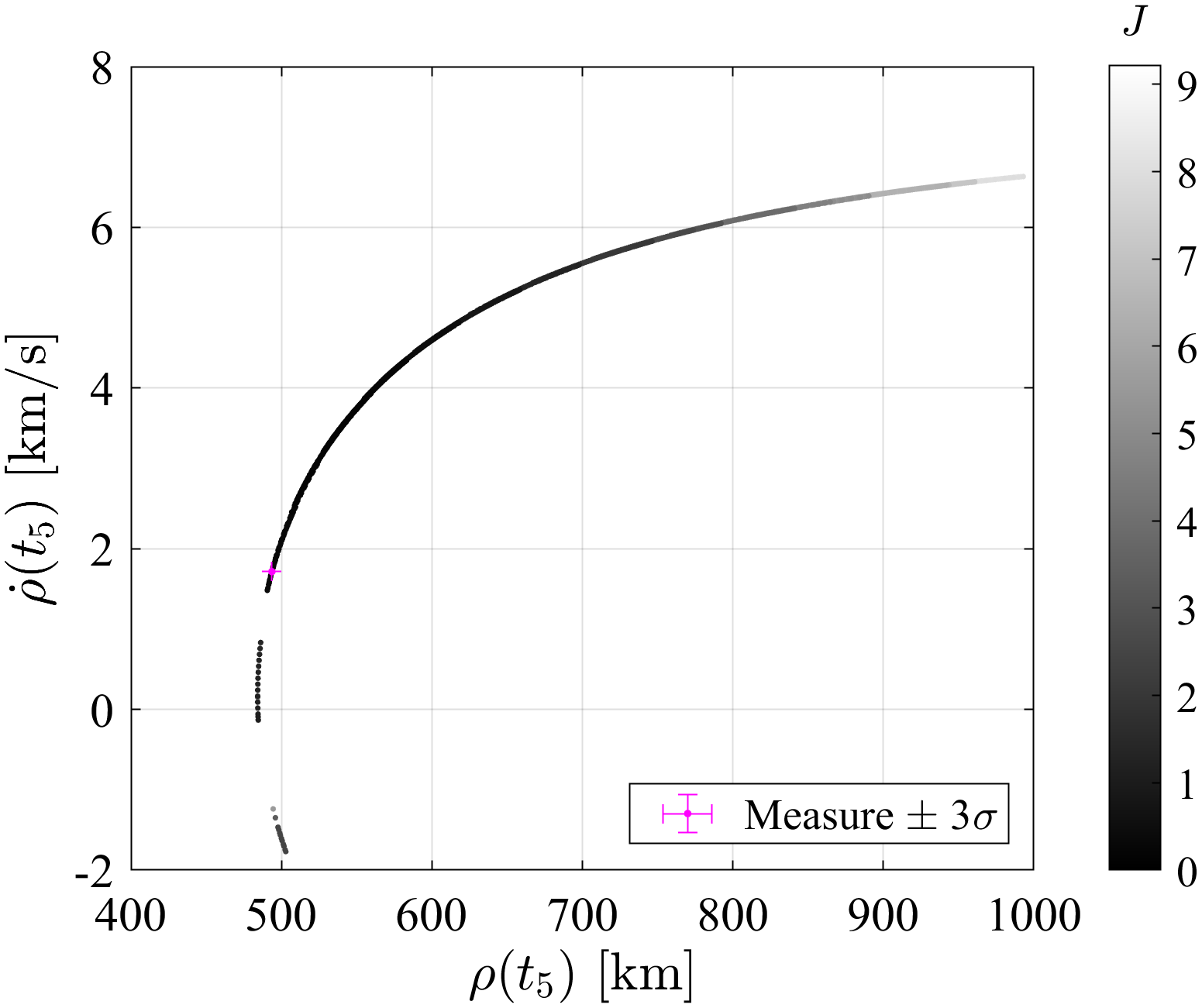}
    \caption{ Object 2.}
    \label{fig:LambertPropag12_TC2}
     \end{subfigure}
     \caption{Projection of points of minimum residual at the next passage, Algorithm 1.}
        \label{fig:ProjectionAlgo1}
\end{figure*}

The states corresponding to the minimum residuals of the retained boxes are propagated forward to the following observation opportunity at $t_5$, which is 11.8 hours later for Object 1 and 12.0 for Object 2. Solutions are projected in the range range-rate plane and shown in Fig. \ref{fig:LambertPropag12_TC1} and \ref{fig:LambertPropag12_TC2}. From these figures, it can be noticed that the new measurement is close to the states corresponding to lower values of the residuals. However, plotting only the minimum values in the domains produces a discontinuous curve that might complicate the automation of the data association process. This drawback is resolved with Algorithm \ref{alg:IOD2}.

\subsection{IOD with accurate angular measurements on one tracklet} % inversion
\begin{figure*}[h]
     \centering
     \begin{subfigure}[b]{0.48\textwidth}
         \centering
       \includegraphics[width=\columnwidth]{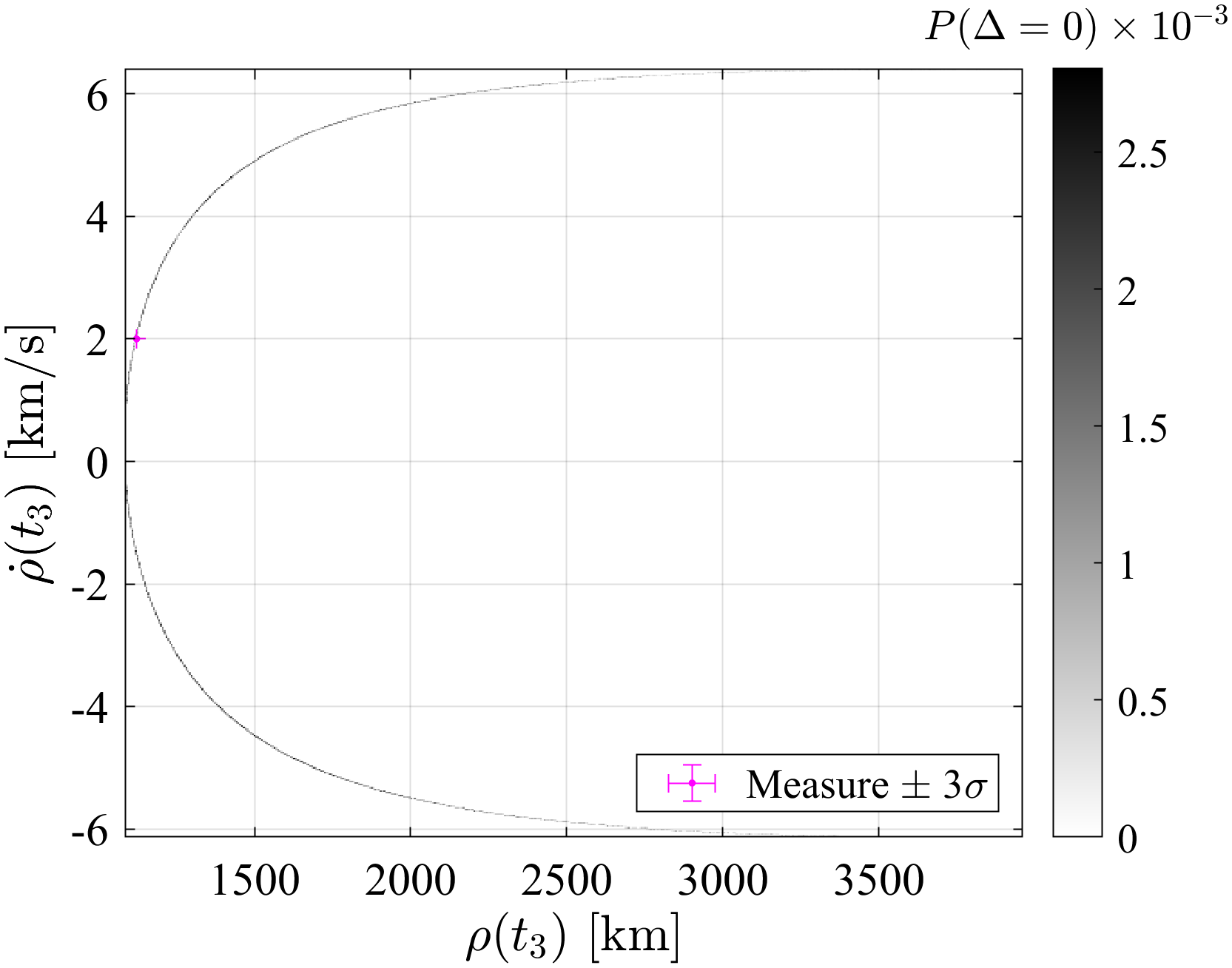}
    \caption{Object 1.}
    \label{fig:InversionHeatmap_TC2}
     \end{subfigure}
     \hfill
     \begin{subfigure}[b]{0.48\textwidth}
         \centering
        \includegraphics[width=\columnwidth]{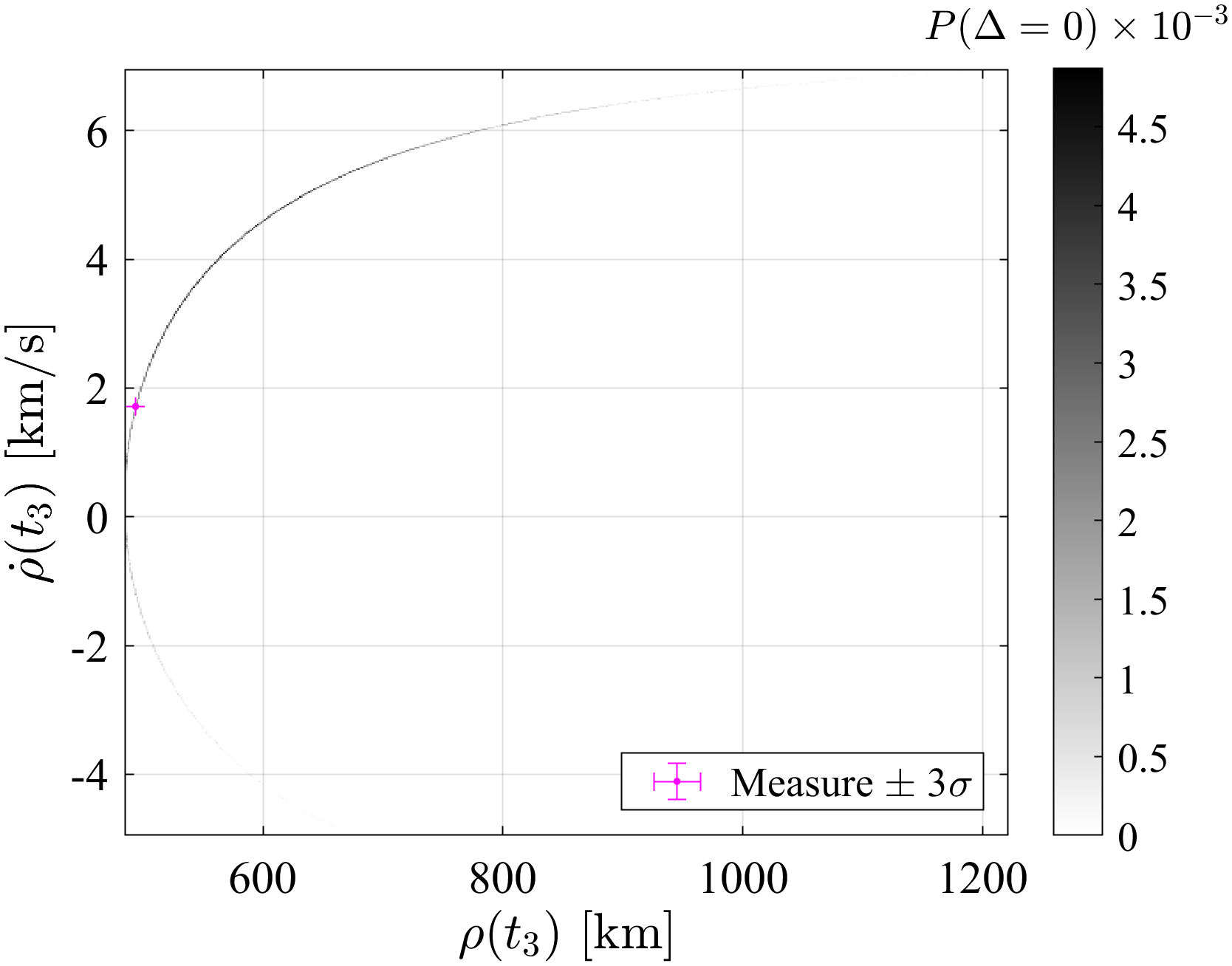}
    \caption{Object 2.}
    \label{fig:InversionHeatmap_TC1}
     \end{subfigure}
     \caption{Projection of points of zero residual at the next passage, Algorithm 2.}
        \label{fig:ProjectionAlgo2}
\end{figure*}

Algorithm \ref{alg:IOD2} is used for accurate measurements on a single tracklet using a significance level $\alpha = 0.01$ on the regressed data. For every value of the observations compatible with their uncertainties,  $({{\rho}}(t_1), \dot{{\rho}}(t_1), {\bar{Az}}(t_1), {\bar{El}}(t_1), {\rho}(t_2), \dot{{\rho}}(t_2))$, the algorithm computes the set of ${\dot{\bar{Az}}}(t_1), \dot{{\bar{El}}}(t_1)$ that result in $\boldsymbol{\Delta} = \boldsymbol{0}$. Similarly to the previous section, this set is propagated forward to the time of the next passage, $t_3$ for data association. Figures \ref{fig:InversionHeatmap_TC2} and \ref{fig:InversionHeatmap_TC1} show heat maps generated by Monte Carlo sampling of the IOD at the time of the subsequent passage. It is evident that, although there is a wide span of possible range and range rate values, the set is defined by a narrow line intersecting the new observations. This is particularly valuable, considering that an analytical $J_2$ model was used for the propagation to simplify the process.

\subsection{IOD with accurate angular measurements on both tracklets} 
\begin{figure*}[!h]
     \centering
     \begin{subfigure}[b]{0.48\textwidth}
         \centering
        \includegraphics[width=\columnwidth]{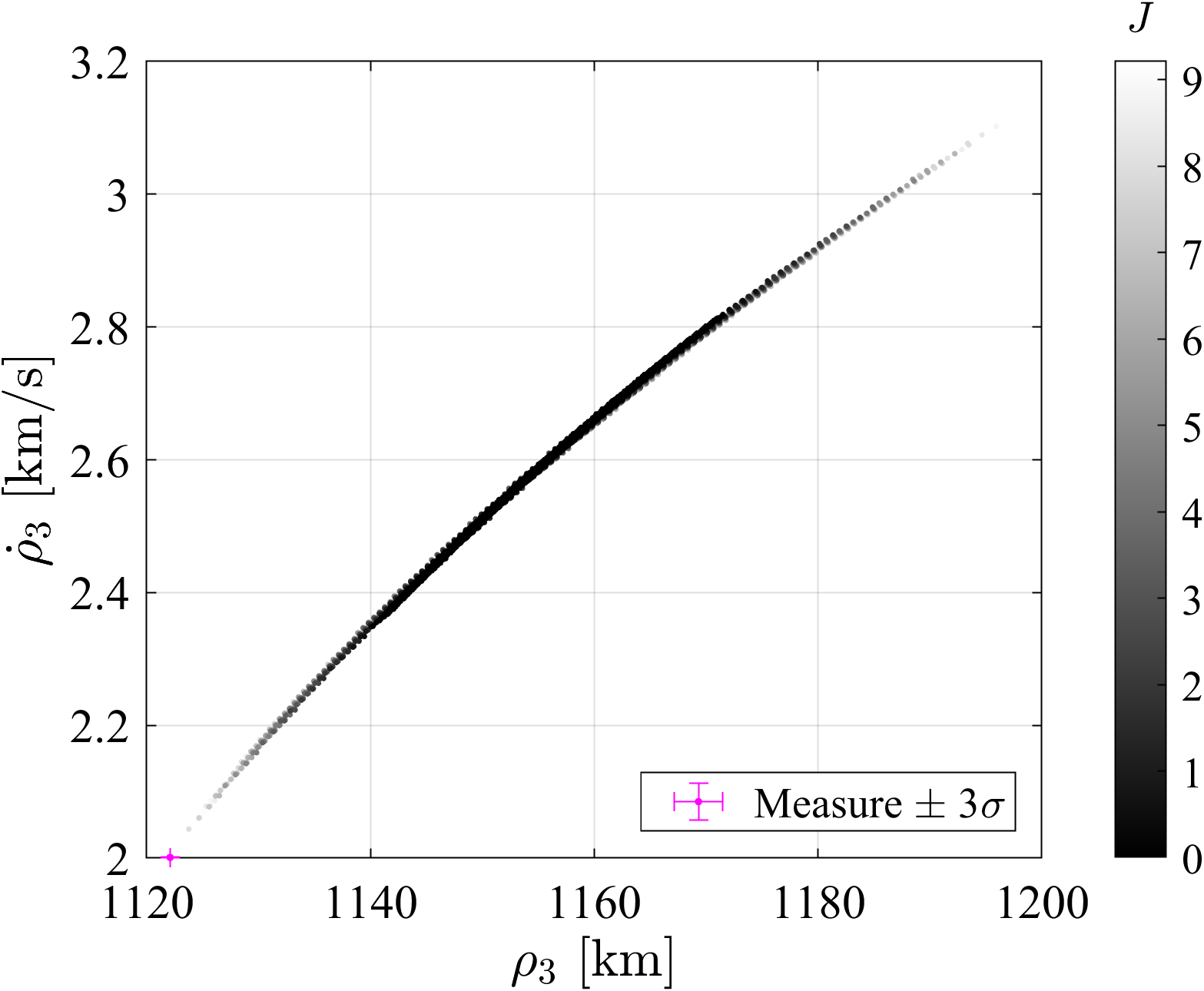}
    \caption{Object 1.}
    \label{fig:InversionHeatmap_TC2_algo3}
     \end{subfigure}
     \hfill
     \begin{subfigure}[b]{0.48\textwidth}
         \centering
        \includegraphics[width=\columnwidth]{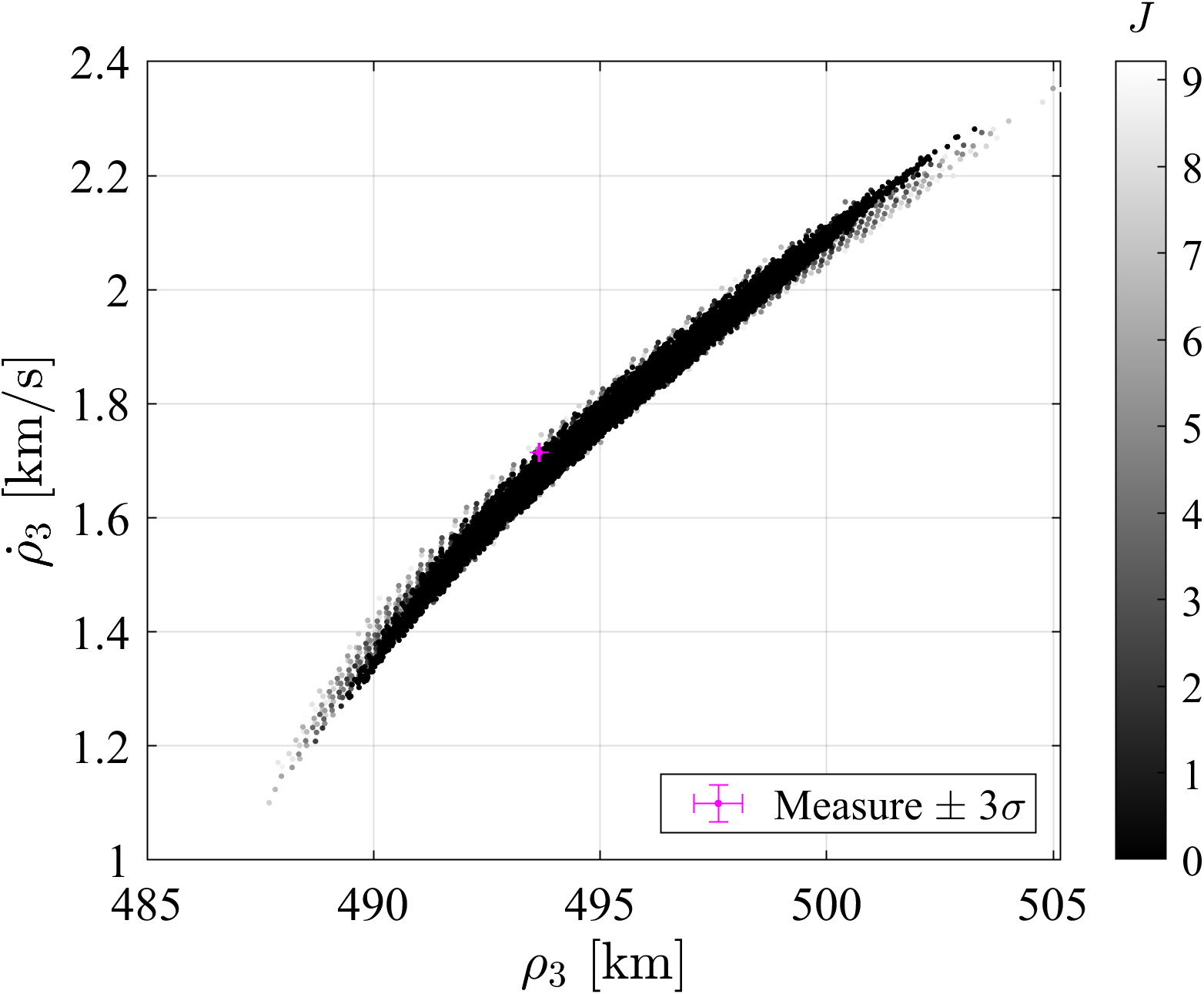}
    \caption{Object 2.}
    \label{fig:InversionHeatmap_TC1_algo3}
     \end{subfigure}
     \caption{Projection of points of minimum residual at the next passage, Algorithm 3.}
\end{figure*}

Accurate measurements for both tracklets allow us to perform a Lambert arc across them. Since these measurements were unavailable for the observed objects on both tracklets, we constructed simulated angular measurements for the second track, stemming from the available state information from the LeoLabs catalog. We then added random noise and assumed the same confidence interval for the available angular measurements. The results of the IOD process in the range-rate plane for the next observation, at $t_3$, are shown in Fig. \ref{fig:InversionHeatmap_TC2_algo3} and \ref{fig:InversionHeatmap_TC1_algo3}. As expected, the distribution is much more constrained than in the previous cases. The plotted points correspond to those for which the normalized squared residual with respect to the range-rate measurements is less than $\chi^2_{2,0.99} \approx 9.21$. However, it is worth noting that for Object 1, the actual measurement does not belong to the computed set. While this method produces a much smaller uncertainty, it is also more sensitive to measurement errors. Therefore, a very good calibration of the angular measurements and high-fidelity dynamical models for propagation is essential.

\begin{figure*}[h]
     \centering
     \begin{subfigure}[b]{0.48\textwidth}
         \centering
       \includegraphics[width=\columnwidth]{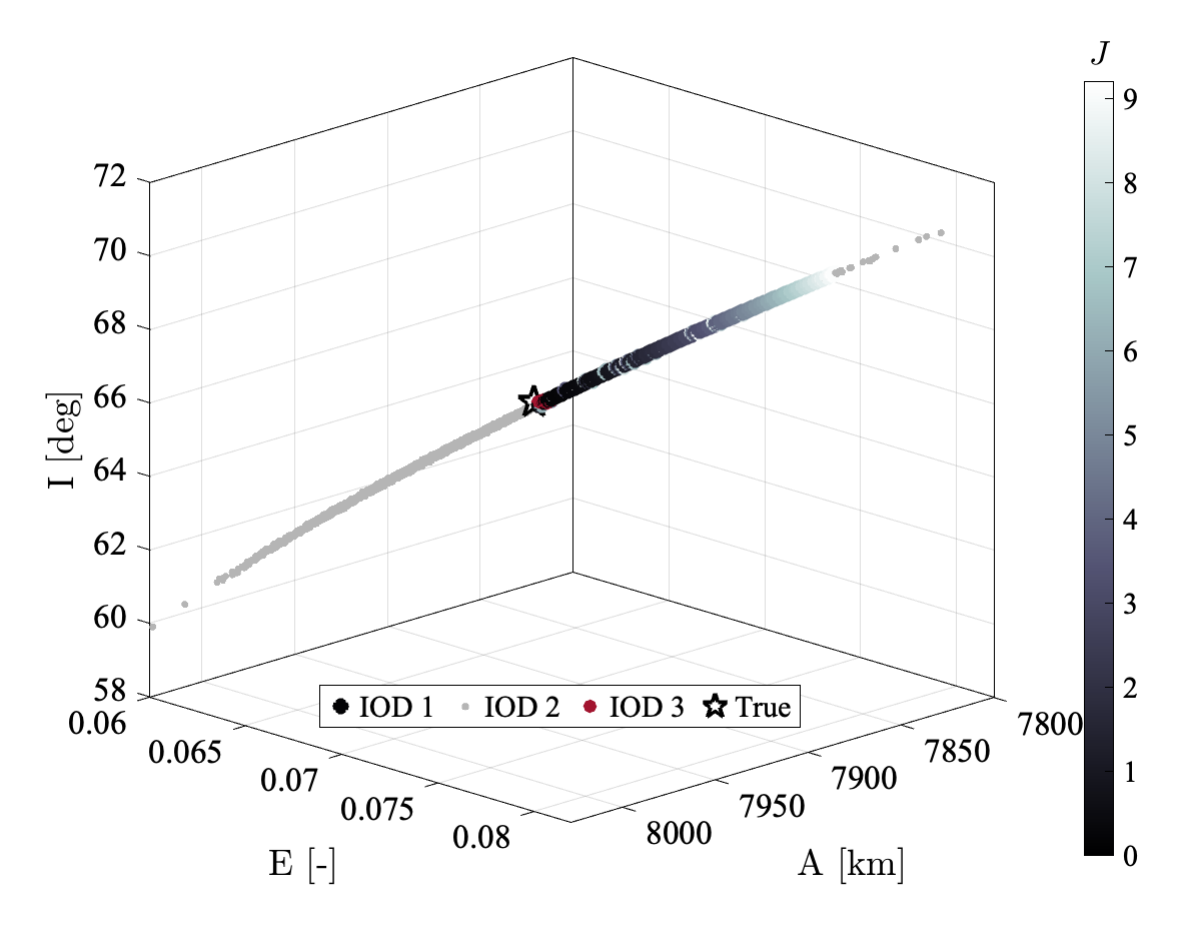}
    \caption{Object 1.}
    \label{fig:3d1}
     \end{subfigure}
     \hfill
     \begin{subfigure}[b]{0.48\textwidth}
         \centering
        \includegraphics[width=\columnwidth]{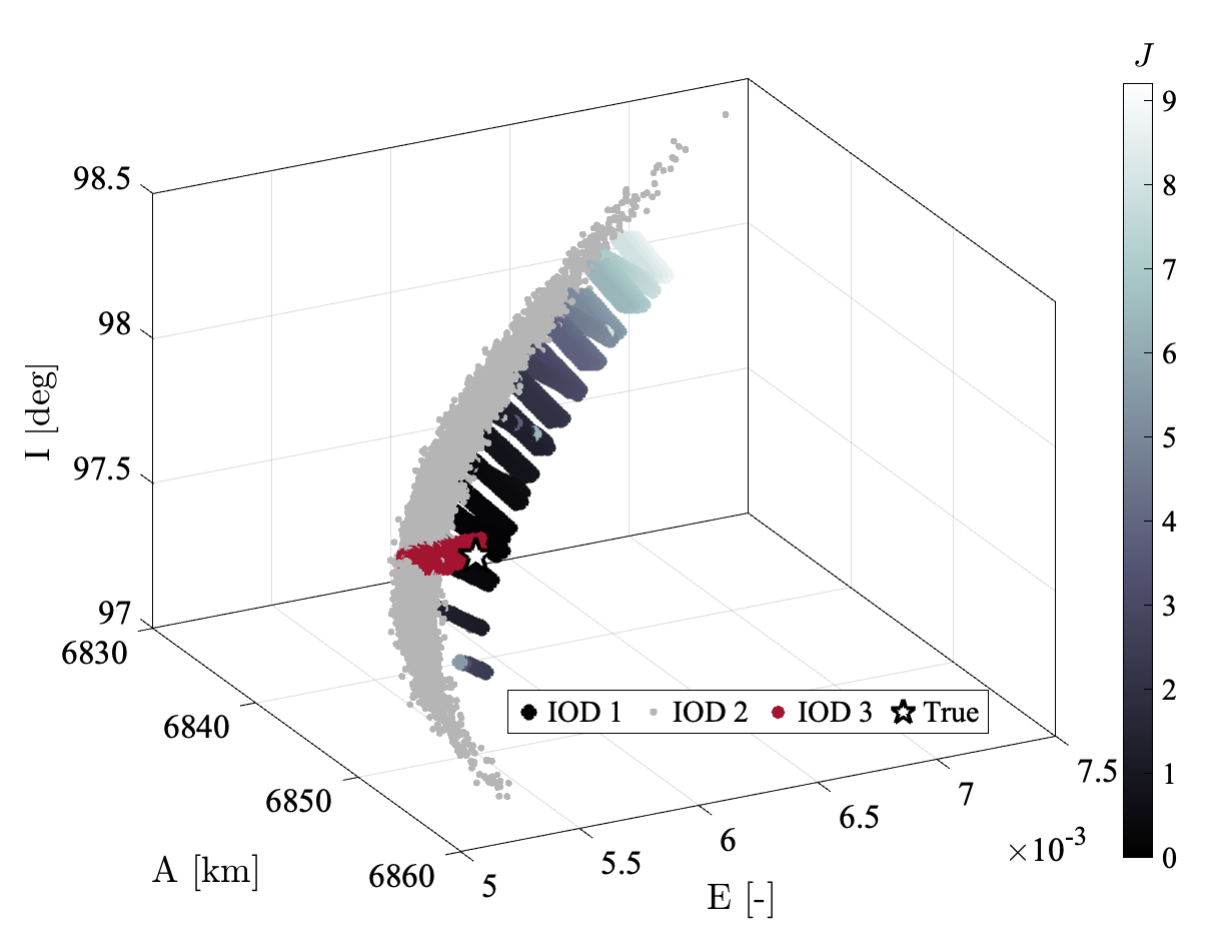}
    \caption{Object 2. }
    \label{fig:3d2}
     \end{subfigure}
     \caption{Uncertainty region for semimajor axis, eccentricity and inclination determined by the three algorithms.}
        \label{fig:3d}
\end{figure*}

To appreciate the three algorithm's roles in exploiting the available data with increasing information going from Algorithm 1 to Algorithm 3 (and hence reducing the final uncertainty on the state), Figures~\ref{fig:3d1} and \ref{fig:3d2} show the uncertainty region in semimajor axis, eccentricity and inclination for the two cases analyzed. The bounds on all orbital elements are reported on \ref{table2}, where it is apparent how the availability of accurate measurements on both tracklets significantly reduced the orbital uncertainty. 

\begin{table*}[h!]
\centering
\begin{tabular}{|c|c|c|c|}
\hline
\multicolumn{4}{|c|}{\textbf{Object 1}} \\ \hline
 & \textbf{Algorithm 1} & \textbf{Algorithm 2} & \textbf{Algorithm 3} \\ \hline
$a$ (km) & [7860.104, 7920.734] & [7860.540, 7994.749] & [7917.445, 7920.857] \\ \hline
$e$ & [0.06985, 0.07754] & [0.06249, 0.07748] & [0.06981, 0.07027] \\ \hline
$i$ (deg) & [65.9, 69.3] & [61.7, 69.3] & [65.9, 66.0] \\ \hline
$\Omega$ (deg) & [146.6, 147.6] & [146.6, 148.9] & [147.6, 147.6] \\ \hline
$\omega$ (deg) & [338.1, 344.7] & [327.3, 344.8] & [338.1, 338.4] \\ \hline
\multicolumn{4}{|c|}{\textbf{Object 2}} \\ \hline
 & \textbf{Algorithm 1} & \textbf{Algorithm 2} & \textbf{Algorithm 3} \\ \hline
$a$ (km) & [6837.281, 6852.650] & [6833.380, 6857.994] & [6847.535,  6848.950] \\ \hline
$e$ & [0.00532, 0.00678] & [0.00505, 0.00726] & [0.00526, 0.00562] \\ \hline
$i$ (deg) & [97.4, 98.3] & [97.1, 98.5] & [97.6, 97.7] \\ \hline
$\Omega$ (deg) & [98.3, 98.5] & [98.3, 98.5] & [98.4, 98.4] \\ \hline
$\omega$ (deg) & [94.1, 126.9] & [88.1, 144.9] & [114.9, 119.7] \\ \hline
\end{tabular}
\caption{Comparison of Parameters for Object 1 and Object 2 across Algorithm 1, Algorithm 2, and Algorithm 3}
\label{table2}
\end{table*}

\section{Conclusions}
\label{Concl}

We developed three algorithms for initial orbit determination that utilize range, range-rate, and angular measurements from LeoLabs radars. These algorithms produce a set of states compatible with the observations represented as a net of Taylor polynomials. 

Our results indicate that data association between very short tracklets lasting a few seconds and separated by two minutes is achievable. IOD solutions can be derived with or without accurate angular measurements when using data from both tracklets. The estimated ranges of orbital elements include or are very close to the actual solution, depending on the algorithm used. When angular measurements are available for both tracklets, the uncertainty reduces significantly. This suggests that the availability of angular measurements for both passages on the co-located radars would be highly beneficial for cataloging purposes. 

In all scenarios, the propagation of uncertain IOD solutions to the next available passage demonstrates that our algorithms can facilitate the data association of uncatalogued objects. This capability presents a valuable opportunity to enhance our understanding of the space debris environment, allowing for the inclusion of new objects in the catalog, particularly smaller and more elusive ones.

However, although not shown here, attempts to use azimuth and elevation rates derived from regressed values proved unsuccessful due to their high sensitivity to measurement errors.

Looking ahead, we plan to test these algorithms on a broader range of objects to validate their effectiveness and expand their applicability in the field of space surveillance. Additionally, we must assess the data association capability in operational scenarios involving multiple uncorrelated tracks that may be compatible with our IOD solutions.
$\,$

$\,$

\bibliography{Ref.bib}

\end{document}